\documentclass{article}
\newtheorem{Def}{Definition}[section]
\newtheorem{Thm}[Def]{Theorem}
\newtheorem{Lem}[Def]{Lemma}
\newtheorem{Cor}[Def]{Corollary}
\newtheorem{Prop}[Def]{Proposition}
\newtheorem{Rem}[Def]{Remark}
\newtheorem{Conj}[Def]{Conjecture}

\usepackage{mathpazo}
\usepackage{amsmath, amssymb}
\usepackage[dvipdfm]{graphicx}
\begin{document}

\title{An Approach to the Hirsch Conjecture} 
 
\date{} 
\author{Yuji Odaka}

\maketitle

\begin{abstract}

W. M. Hirsch proposed a beautiful conjecture on diameters of convex polyhedra, 
which is still unsolved for about 
50 years. I suggest a new method of argument from the viewpoint of deformation and moduli of 
polytopes. As a consequence, for example, if there are at least $3$ \textit{disjoint geodisics} 
for all Dantzig figures, as in the $3$ dimensional case, the conjecture follows.
 
\end{abstract}

\section{Introduction}
The graph-theoritic study of polytopes has a long history, as the objects are very natural and so interesting.  
The aim of this short paper is to suggest an approach to one of the fundamental, important, and long-standing 
conjecture in that area. 
From now on, 
we are always concerned about the graph-theoritic structure of the polytopes, and call diameters etc in that sense. 

In 1957, W. M. Hirsch proposed the Hirsch conjecture and the $d$-step conjecture. The Hirsch conjecture states 
that the maximum of possible diameters of $d$-polyhedra with $n$ facets, where $d\geq 2$ denotes the dimension,  
which is denoted by $\Delta (d, n)$ is less than or equal to $n-d$, and as a special case of $n=2d$, the $d$-step 
conjecture states $\Delta (d, 2d)=d$. Though these are proved to be 
equivalent, they are unsolved. We have to pay attention that there are $2$ versions to each 
conjecture, corresponding to whether we admit unbounded polyhedron or not. Only bounded polyhedron will 
be  
called \textit{polytope} in this paper. We will consider only that   
case. 

For the evaluation of $\Delta (d, n)$, 
Adler [2] proved the left hand side of following inequality, and Kalai-Kleitman [6] recently 
proved the right hand side. 

\begin{equation}
(n-d)-\frac{n-d}{[5d/4]}-1\leq \Delta (d, n) \leq n^{\log{d}+2}
\end{equation}
 
However, they do not tell even whether  $\Delta (d, 2d)$ grows in polynomial order. 

The most simple equivalent form of this conjecture is the \textit{non-revisiting path 
conjecture}, whose equivalence to the original is proved in [1]. It conjectures the 
following: any two vertices of a simple polytope $P$ can be joined by a path that does not revisit any facet 
of $P$. Actually in [1], more is proved. The $d$ dimensional Dantzig figure, defined in [1], 
means the 
triplet $(P,x,y)$ where $P$ is a $d$ dimensional simple convex polytope, with just $2d$ facets, 
half of which includes $x$, and the others includes $y$. In [1], it is proved that if for all the simple Dantzig 
figure $(P,x,y)$, the distance between $x$ and $y$ is not more than its dimension, 
the Hirsch 
conjecture holds. When there are no problems, we sometimes omit $(x, y)$
. Furthermore, in [1], the bounded version of Hirsch conjecture has 
already proved for the dimension less than $6$, and there is counterexample proposed for the unbounded case. 
For more details about the conjecture, I recommend a great summary [4]. 

In this paper, we call a polytope (resp.Dantzig figure) \textit{Hirsch-polytope} (resp.  \textit{Dantzig-Hirsch polytope}), 
if it is a convex polytope that satisfies the conjecture. And we say 
the polytope is $(n, d)$-type, when it is $d$ dimensional with $n$ facets. 

\section{The Fundamental Deformation}
The polytopes are in $d$ dimensional Euclidean space $\mathbb{R}^d$, 
and each of the 
facets are the parts of hyperplanes in $\mathbb{R}^d$. Note $F_{i}[P](1\leq i\leq n)$ the facets of the 
particular polytope$P$, $\pi_{i}[P]$ the hyperplane corresponding to the facet $F_{i}[P]$. We often omit it  
as $F_{i}$. 

We introduce the basic concept and the key of this paper, \textit{fundamental deformation}. 

\begin{Def}
Think of the continuous move of a particular hyperplane $\pi_{i}$, which corresponds to a facet $F_{i}$ of P. 
Strictly speaking, this means the family of polytopes $P(t)(0\leq t\leq 1)$, whose facets are fixed     
$\pi_{j}[P]$s 
with $i\not=j$   
and a hyperplane $\pi_{i}(t)(0\leq t\leq 1)$ with $\pi_{i}(0)=\pi_{i}$, which is continuous family. If the 
move 
satisfies the following, we call it \textit{fundamental deformation}. We will often omit as FD. 
\begin{enumerate}
\item{$P(t)(0\leq t\leq 1)$ are simple polytopes except for one $t$, which we denote $t_0(<1)$. }
\item{there is only one vertex $v$ of $P$ which is an intersection of some $\pi_{j}$s with $i\not=j$, passed 
through by 
the move $\pi_{i}(t)$. Actually, $\pi_{i}(t_0)$ passes the point $v$, and it is not on $P(t)$ with $t>t_0$
. }
\item{$P(t)\subset  P(s)$ for $t>s$} 

\end{enumerate}
From the conditions, $P(t)(0\leq t<t_0)$ have all the same combinatoric type of polytope, and so all $P(t)$ 
 with $t_0\leq 
t<1$ are. We call \textit{deformation} the iteration of fundamental deformations. Sometimes, we call $P(1)$ itself a 
fundamental deformation. Just one edge will vanish in the process of fundamental deformation. We call it 
vanishing edge.
\end{Def}

\begin{Rem}
From now on, we assume $d\geq 3$ so that the fundamental deformation of (n, d)-type polytope has (n, d)-type too. 
It is easy to see that it is true iff $d\geq 3$. 
\end{Rem}
\begin{Def}
The \textit{FD of Dantzig figure} is a bit different notion. Dantzig figure is the \textit{triplet}
 so we define it as the FD as in the above sense whose vanishing edge does not includes 
$\{x,y\}$ as its vertices. In the part of following argument where we will treat only Dantzig figures, we 
omit FD of Dantzig figures just as FD.
\end{Def}

We will see the detail of the structure of the fundamental deformation. Notation remains and we omit the proof but it is easy: 
\begin{Prop}
Let $P(t)$ be the FD. 
As $P(0)$ is simple, we can assume that ${v}$ is the intersection of $\pi_j(2\leq j\leq d+1)$ and $i=1$. 
Permute the order of $\pi_j(j\geq 2)$ if necessary, so that $v$ and $w$, which is the intersection of 
$\pi_1(t)$ and $\pi_i(3\leq i\leq d+1)$ are only points which will vanish 
in the process of FD. 

Now we can describe the basic property of FD as follows. 
\begin{itemize}
\item{For $0\leq t<t_0$, the segment $vw$ is an edge of the polytope, that is the intersection of 
$\pi_j(3\leq j\leq d+1)$. }
\item{For $t_0<t\leq 1$, the intersection of $\pi_j(j=1, 2)$ and the $P(t)$ is $d-2$-simplex. }
\end{itemize}
\end{Prop}

What is the most important idea in this paper is the following:
\begin{Cor}
If $x, y$ are distinct vertices on the polytope and $p:[0, l]\rightarrow P(=P(0))$ be the edge-following 
path (i.e. 
 $p(s)$ is vertex$\Leftrightarrow$ $s$ is integer and ${p(x)|m\leq x\leq m+1}$ for integer $m$ is a segment 
of $P$. $l$ will be called the 
\textit{length} of $p$. )If the following condition \textit{does not} holds, there is  an edge-following path 
$p'$ 
of 
$P(1)$ with the same or the fewer length with $p$. 

\begin{itemize}
\item{$p$ passes through just one of $\{v,w\}$.}
\end{itemize}
In the proof, we construct $p'$. We will call it the fundamental deformation of $p$.   
\end{Cor}
\textbf{proof of the Corollary:} 
We can assume that if $p$ passes both $v$ and $w$, it passes the edge $vw$. 
Then, there are two cases possible for 
the edge $p(n)p(n+1)$. 
\begin{itemize}
\item{If it does not touch the segment $vw$, we does not change it. }
\item{If it is $vw$ itself, we can take alternative segment without modifying $p(n-1)p(n)$ and $p(n)p(n+1)$
(strictly speaking, the part of them) from the $d-2$simplex which will be constructed between $\pi_{1}$ and 
$\pi_{2}$ in $P(1)$. We can do it since the diameter of the simplex is $1$. }
\end{itemize}
$\blacksquare$
\section{Our Program}

In this section, we will explain the fundamental principle. 

Let $G(n, d)$ be the set of combinatoric-type of all the polytopes of $(n, d)$-type. Regard this 
as an oriented graph, with vertices set itself, and arrows set naturally corresponding to fundamental 
deformations. 

From the previous corollary, we see the basic and most important property of this graph: 
\begin{Def}
We introduce the concept of \textbf{good} fundamental deformation(resp.of Dantzig figure) here.
It means the FD satisfying the following condition:
\begin{itemize}
\item{For any pair of points $\{x,y\}$ (resp.For $\{x,y\}$ of the original Dantzig figure $(P,x,y)$), 
there is a path connecting them with minimum length either passes the vanishing edge or does not even touch 
it.} 
The deformation which is an iteration of good fundamental deformation is said to be good.
\end{itemize}
\end{Def}
\begin{Lem}
The good deformation of Hirsch (resp.Dantzig-Hirsch polytope) is also Hirsch (resp.Dantzig-Hirsch polytope).
\end{Lem}
\textbf{proof:} Obvious from the previous Corollary. $\blacksquare$ 

We restrict our attention to the $d$ dimensional 
\textit{Dantzig figures}. Let $D(d)$ be the set of possible combinatoric types of $d$ dimensional Dantzig 
figures. And regard it as graph as we did for $G(n, d)$. We call these, 
\textit{moduli-graph}s. 

\begin{Prop}
$D(d)$ is strongly connected graph; for any two $d$ dimensional Dantzig figures $D_1$ and $D_2$
, there is  other $d$ dimensional Dantzig figures $D'_1$ and $D'_2$ such that the latter is a deformation of 
the former one. 
\end{Prop}
\textbf{proof:}
\begin{Lem}
For all $d$ dimensional Dantzig figures $(D, x, y)$, there is a Affine transformation $f$ such that:

\begin{itemize}
\item{$f(x)=(0, \dots, 0)$ and $f(y)=(0, \dots, 0, 1)$}
\item{$f(D)\subset \{(x_1, x_2, \dots, x_n)\mid 0\leq x_n \leq 1\}$}
\item{$f(D)\cap \{(x_1, x_2, \dots, x_n)\mid  x_n=0$ or $ 1\}={x, y}$}
\end{itemize}
\end{Lem}
\textbf{proof of the Lemma:} As $D$ can be written as an intersection of two cones $C_x$ and $C_y$, which have 
$x$ and $y$ respectively as 
vertices. At first we linear-transform $D$ to make $C_x$ identical to $\{(x_1, \dots, x_n)\mid x_i\geq 0\}+x$. 
Then, the boundedness of $D$ implies just $C_y\cap\{y+(x_1, \dots, x_n)\mid x_i\geq 0\}={(0, \dots, 0)}$. 

From 
this argument with coordinate, we can see now $C_x-x$ and $C_y-y$ (parallel transformation of the original 
cones) have only the origin in common. Therefore from the famous \textit{separation theorem} (cf. 
[7]), there is  a hyperplane $\pi$ which have origin as only common point with $C_x-x 
\cup C_y-y$, 
separating them. Then, the parallel hyperplane of $\pi$ which pass through $x$ and $y$ put $D$ between them. 

Finally, as we can see easily, there is  an Affine transformation which moves $x$ and $y$ to $(0, \dots, 0)$ and 
$(0, \dots, 0, 1)$ respectively conserving the condition that \textit{parallel planes which pass through $x$ and 
$y$ as only common point with $D$ respectively, put $D$ between them. }   
$\blacksquare$

\textbf{proof of the Proposition (continued):}

Take the lemma's Affine transformation $f_1$ and $f_2$ of $D_1$ and $D_2$ respectively so that they have same 
positions of $x$ and $y$. It is easy to see that 
we can take them so that $f_1(D_1)\supseteq  f_2(D_2)$ and we replace $D_i$ by $f_i(D_i)$. Then, let us assign 
each facets of $D_1$ which pass 
through $x(resp. y)$, a facet of $D_2$ which passes through $x(resp. y)$ bijectively so that orientation-
preserving homeomorphism from the neighbor of $x$(resp.$y$) in $D_1$ to that in $D_2$ which send the parts of 
facets of $D_1$ to corresponding facets of $D_2$. Just one by one, rotate the facets of $D_1$ one way around 
the intersection of $2$ corresponding facets which has codimension$2$, we complete the proof of the 
proposition here.  
$\blacksquare$

Now, we can states our program to prove the Hirsch conjecture. First in the rough form. We may call this 
\textit{Deformation  
program}. 
\begin{quotation}
Think of the subgraph $E(d)$ of $D(d)$ with the same vertices set and the edge sets consist of the original 
one 
corresponding to the \textbf{good} fundamental deformations. Find \textit{many} Hirsch-polytopes, and show that every vertex is 
a destination of oriented path in $E(d)$ from one 
of them.   
\end{quotation}  
\section{3-geodisic conjecture}
Previous section is just an rough idea. 
In this section, we will introduce some new conjectures which follows those basic ideas. 

The basic conjectures are the following:
\begin{Conj}[Dantzig conjecture]
$E(d)$ is connected in strong sense. 
i. e. For any vertex $x, y$ of the graph $E(d)$, there are (oriented) path from $x$ to $y$.  
\end{Conj}
\begin{Conj}[Strong Dantzig conjecture] 
$E(d)=D(d)$ i.e. all the fundamental deformation of $d$ dimensional Dantzig figure is good. 
\end{Conj}
From what we said in the previous section, $D(d)$ is connected in the strong sense. Therefore, the Strong Dantzig conjecture
derives the Dantzig conjecture. 
\begin{Def}
If $G$ is the graph (oriented (resp. unoriented)), we denote $m(G:x,y)$ the maximum number of point-distinct 
geodisics(oriented(resp. unoriented)). Point-distinct means any $2$ of them only have common 
point $x, y$
, and geodisic means the path with minimum length between the particular $2$ points. 
\end{Def}
Then the following holds. 
\begin{Thm}
If for all the $d$ dimensional Dantzig figures, $m(P:x,y)\geq 3$ holds, then the Strong Dantzig conjecture 
is true for that dimension. 
\end{Thm}
\textbf{proof:}
Let $(P,x,y)$ be the fixed $d$ dimensional Dantzig figure. 
It is sufficient to prove that for any edge $e$ of $P$ which does not touch $x$ nor $y$, 
there is a geodisic between $x$ and $y$, which satisfies one of the following: 
\begin{itemize}
\item{$p$ does not touch $e$. }
\item{$p$ includes $e$ (not just touching)}
\end{itemize}  

The edges which satisfys one of the above will also be called \textit{good} and not good edge which does not 
touch $x$ nor $y$, will be called \textit{bad}. Assume there is  bad edge $e_0$, and let us $p_1$, $p_2$ be any of 
the geodisics between $x$ and $y$. By definition, $e_0$ should have one vertex on $p_1$ and the other on the 
$p_2$. So if there is  $3$ geodisics between $x$ and $y$, it contradicts. 
$\blacksquare$

\begin{Conj}[3-geodisic conjecture]
For $d$ dimensional Dantzig figure $(P,x,y)$, $m(P:x,y)\geq 3$. 
\end{Conj}
The theorem states, this is stronger than the Strong Dantzig conjecture. 

\begin{Rem}
Is this true? 
The first nontrivial case is $d=3$. In this case, we can easily confirm all the conjectures hold.
We have not checked even for $d=4$. If we take off the geodisic condition, it is well-known that, 
the number is at least $d$ ,so more than or equal to $3$. The geodisic condition is very strong, but to be a 
Dantzig figure is a restrictive 
condition...
\end{Rem}
\textbf{Acknowledgement}

I would like to thank my friends 
 K.Yamashita, K.Irie and S.Ohkawa for some comments on the manuscript of this paper.

  Yuji Odaka: RIMS, Kyoto University, Oiwake-Cho, Kitashirakawa, Sakyo-ku, Kyoto 606-8502, Japan

\textit{Email address:} \textbf{ yodaka@kurims.kyoto-u.ac.jp}

\end{document}